\documentclass{amsart}
\usepackage{mathrsfs}
\usepackage{amsfonts}
\usepackage{amssymb}
\usepackage[papersize={7.3in,10in},textwidth=14cm,textheight=20.6cm,centering]{geometry}
\usepackage{enumerate}

\usepackage[colorlinks=true,citecolor=blue,linkcolor=blue]{hyperref}
\hypersetup{
pdfstartpage=1,
pdfstartview=FitH}

\allowdisplaybreaks

\newtheorem*{theorem*}{Theorem}
\newtheorem{theorem}{Theorem}
\newtheorem{lemma}{Lemma}

\newtheorem{proposition}{Proposition}

\newtheorem{conjecture}{Conjecture}

\newtheorem*{claim*}{Claim}

\theoremstyle{remark}

\newcommand{\ud}{\mathrm{d}}
\newcommand{\sym}{\mathrm{sym}}

\DeclareMathOperator*{\Res}{Res}
\DeclareMathOperator{\Mod}{mod}
\DeclareMathOperator{\Frob}{Frob}
\DeclareMathOperator{\tr}{tr}
\renewcommand{\bmod}[1]{\,(\Mod{ #1})}

\begin{document}

\title{Sign changes of Kloosterman sums with almost prime moduli}
\thanks{The work is partially supported by China Scholarship Council and N.S.F. (No. 11171265) of P. R. China.}

\author{Ping Xi}

\address{School of Mathematics and Statistics, Xi'an Jiaotong University, Xi'an 710049, P.R. China}
\address{EPFL/SB/MATHGEOM/TAN, Station 8, 1015-Lausanne, Switzerland}

\email{pingxi.cn@gmail.com,\ ping.xi@epfl.ch}


\date{\today}


\begin{abstract}We prove that the Kloosterman sum $S(1,1;c)$ changes sign infinitely often as $c$ runs over squarefree moduli with at most 10 prime factors, which improves the previous results of E. Fouvry and Ph. Michel, J. Sivak-Fischler and K. Matom\"{a}ki, replacing 10 by 23, 18 and 15, respectively. The method combines the Selberg sieve, equidistribution of Kloosterman sums and spectral theory of automorphic forms.
\end{abstract}

\maketitle

\tableofcontents

\section{Introduction}\label{sec:1}
In this paper, we are interested in the sign changes of
Kloosterman sums,
defined by 
\[S(m,n;c)=\sideset{}{^*}\sum_{a\bmod c}e\bigg(\frac{ma+n\overline{a}}{c}\bigg)\]for each positive integer $c$ and integers $m,n$, where $a\overline{a}\equiv1\bmod c.$
There is much literature investigating the Kloosterman sums
because of their profound applications
in analytic number theory and automorphic forms
as well as their own mysterious nature.

A well-known estimate for individual Kloosterman sums due to
A. Weil \cite{We}
asserts that
\begin{align}\label{eq:1}|S(m,n;p)|\leqslant2p^\frac{1}{2}\end{align}for each prime $p$ with $(m,n,p)=1.$
More generally, one has
\[|S(m,n;c)|\leqslant c^\frac{1}{2}(m,n,c)^\frac{1}{2}\tau(c),\]
where $\tau(c)$ is the divisor function;
in fact, T. Estermann \cite{Es} showed the
slightly stronger estimate
\begin{align}\label{eq:2}|S(m,n;c)|\leqslant c^\frac{1}{2}(m,n,c)^\frac{1}{2}2^{\omega(c)}\end{align}
for $32\nmid c$,
where $\omega(c)$ denotes the number of distinct prime factors of $c$. 

Kloosterman sums have
long been basic tools in the analytic theory of
automorphic forms;
for example,
they appear
in the Petersson trace formula for the average of
products of Fourier coefficients of holomorphic modular forms.
In return, the theory of automorphic forms can be used
to study Kloostermans sums.
The precise link was first established by Kuznetsov \cite{Ku},
who, by
means of his trace formula, made progress on a conjecture of
Linnik and Selberg that
\begin{align}\label{eq:3}\sum_{c\leqslant x}\frac{1}{c}S(m,n;c)
  = O_{m,n,\varepsilon}(x^\varepsilon)
\end{align}
for any
$\varepsilon>0$.
Kuznetsov proved that \eqref{eq:3}
is valid for any $\varepsilon>1/6$,
while applying (\ref{eq:2}) to each summand
gives \eqref{eq:3}
only for $\varepsilon > 1/2$.

One might expect that Kuznetsov's estimate is mainly due to the
oscillations of Kloosterman sums as $c$ varies amongst the
consecutive integers,
but one had to wait for the work of
Ph. Michel \cite{Mi}, who was able to confirm this phenomenon by
proving that there must be a positive portion of prime pairs $(p,q)$
such that $|S(1,1;pq)|\geqslant 0.64\sqrt{pq}$. Hence it is
natural to investigate the sign changes of Kloosterman sums when
$c$ varies over thinner set, for instance, the primes.

As an analog of the celebrated Sato-Tate conjecture for elliptic curves, N.M. Katz \cite{Ka1} formulated a conjecture for the equidistribution of the Kloosterman sum angle $\theta_p(a)$, which is defined as 
\[S(a,1;p)=2p^\frac{1}{2}\cos\theta_p(a)\]
by means of (\ref{eq:1}).
\begin{conjecture}[N.M. Katz]\label{conj:1}For
  any
  $f\in\mathcal{C}([0,\pi])$ and 
  nonzero integer $a$, we have
\[\lim_{x\rightarrow+\infty}\frac{1}{\pi(x)}\sum_{p\leqslant x}f(\theta_p(a))=\frac{2}{\pi}\int_0^\pi f(\theta)\sin^2\theta\ud\theta.\]\end{conjecture}

This conjecture predicts that for such an $a$
the angles $\theta_p(a)$ equidistribute with respect to the Sato-Tate measure 
$$\mu_{\rm ST}=\frac{2}{\pi}\sin^2\theta\ud\theta$$
as $p$ runs over all the primes;
it would then follow immediately that $S(1,1;p)$ changes
sign infinitely often as $p$ varies.

There are many facts that support Conjecture \ref{conj:1}. For
instance, Katz himself \cite{Ka2} proved that $\{\theta_p(a) : a \in \mathbb{F}_p^\times\}$
equidistributes with respect to the Sato-Tate measure $\mu_{\rm ST}$ as $p$ tends
to infinity;
we will come back to this issue in the next section.
It is also known
that $S(1,1;c)$ change signs infinitely often as $c$ runs over
positive squarefree integers with at most 23 prime factors,
or more precisely that
\begin{align*}|\{X<c\leqslant2X:S(1,1;c)\gtrless0,\mu^2(c)=1,\omega(c)\leqslant23\}|\gg\frac{X}{\log X}.\end{align*}
This
was proved by E. Fouvry and Ph. Michel \cite{FM1,FM2}
by a pioneering combination and application of the Selberg sieve, spectral theory of automorphic forms and $\ell$-adic cohomology.
The subsequent improvements are due to J. Sivak-Fischler \cite{SF1,SF2} and K. Matom\"{a}ki \cite{Ma}, who reduced 23 to 18 and 15, respectively. 

In this paper, we shall present a further improvement on the problem of sign changes. We would also use the Selberg sieve, but with a modification, inspired by an old idea of Selberg \cite{Se} towards the Twin Prime Conjecture. This will be explained in the next section.

The main theorem can be stated as follows.
\begin{theorem}\label{thm:1}
  There exists an absolute constant $c_0>0$
  such that
  for sufficiently large $X>0$,
\begin{align*}|\{X<c\leqslant2X:S(1,1;c)\gtrless0,\mu^2(c)=1,\omega(c)\leqslant10\}|\geqslant c_0\frac{X}{\log X}.\end{align*}
\end{theorem}

\noindent {\bf Notation.} Throughout this paper, $p$ is reserved for a prime number; we write $e(z)=e^{2\pi iz}$; $\mu,\varphi$ denote the M\"{o}bius and Euler functions, respectively, $\tau$ denotes the divisor function, and $\omega(n)$ denotes the number of distinct prime divisors of $n$. Moreover, $(a,b)$ and $[a,b]$ denote the g.c.d. and l.c.m. of $a,b$, respectively. Given $X\geqslant2$, we set $\mathcal{L}=\log X$. We use $|\cdot|$ to denote the cardinality of a set or the absolute value of a  number.  We adopt the notation $(\sigma)$ to denote the usual contour integral over the line $\sigma+it,t\in\mathbb{R}.$ We use $A$ to denote a sufficiently large positive number and $\varepsilon$ a sufficiently small positive number, which can be different at each occurrence. 

\bigskip

\noindent \textbf{Acknowledgement.} The present work in this paper will be part of my PhD thesis. I am grateful to Professor Philippe Michel for his kind supervision and suggesting this problem to me. His valuable suggestions and comments should be greatly acknowledged. I also thank Paul Nelson for his helpful comments on sieve theory and everything else. The idea in this paper was partially inspired by a talk of Professor Kai-Man Tsang during a conference celebrating the 25 years Number Theory Seminar at ETH Z\"urich in June 2013. I would like to thank Professor Tsang and the organizers of the conference. Sincere thanks are also due to Professor Yuan Yi for her constant help and encouragement. The numerical computations in this paper are based on the Mathematica codes of Kaisa Matom\"{a}ki, and I thank her for sharing the codes on her homepage. I am also grateful to the referee for his/her detailed comments and suggestions, which have greatly improved the exposition of the paper.

\bigskip

\section{Outline of the proof}\label{sec:2}

We prove Theorem \ref{thm:1} by
applying the Selberg sieve. 
Let
$\lambda=(\lambda_d)$ be the Selberg sieve weight given by
\begin{align*}\begin{cases}\lambda_1=1,\\
|\lambda_d|\leqslant1,\\
\lambda_d=0,\ \ \text{if $d>\sqrt{D}$ or $\mu(d)=0$.}\end{cases}\end{align*}
Here $\sqrt{D}=X^\gamma\exp(-\sqrt{\mathcal{L}})$ for some
$\gamma\leqslant\frac{1}{4}$
to be optimized later, and
\begin{align}\label{eq:4}\lambda_d=\mu(d)\bigg(\frac{\log(\sqrt{D}/d)}{\log \sqrt{D}}\bigg)^k\end{align}for $1\leqslant d\leqslant  \sqrt{D}$ and $k$ a positive integer to be specialized later.

Let $g(x)$ be a fixed smooth function supported in $[1,2]$,
and its Mellin transform is defined as
\begin{align*}\widetilde{g}(s)=\int_0^{+\infty}g(x)x^{s-1}\ud x.\end{align*}
Integrating by parts, we have
\begin{align*}\widetilde{g}(s)\ll(|s|+1)^{-A}\end{align*}
for any $A\geqslant0.$

Our starting point is the following sum
\begin{align}\label{eq:5}H^\pm(X)=\sum_ng\left(\frac{n}{X}\right)\frac{|S(1,1;n)|\pm S(1,1;n)}{\sqrt{n}}\mu^2(n)\bigg(\rho-\bigg(\frac{k}{2}\bigg)^{\omega(n)}\bigg)\bigg(\sum_{d|n}\lambda_d\bigg)^2,\end{align}
where $\rho$
is a parameter (depending upon $k$) to be chosen later. Our basic strategy is to show that there exists some pair $(k,\rho)$ with $k\geqslant3$ and $\rho>1$, such that
\[H^\pm(X)>0\]
for $X$ large enough;
it then follows from the definition
that there exists $n\in(X,2X]$ with
\[\omega(n)\leqslant\bigg[\frac{\log\rho}{\log(k/2)}\bigg]\]
for which $S(1,1;n)\gtrless0$. More precisely, one can obtain
a lower bound for the number of such $n$ by applying
H\"{o}lder's inequality appropriately;
this will establish Theorem \ref{thm:1}.

From (\ref{eq:5}), we have
\begin{align*}H^\pm(X)\geqslant\rho H_1(X)-2H_2(X)\pm\rho H_3(X),\end{align*}
where
\begin{align*}H_1(X)&=\sum_ng\left(\frac{n}{X}\right)\frac{|S(1,1;n)|}{\sqrt{n}}\mu^2(n)\bigg(\sum_{d|n}\lambda_d\bigg)^2,\\
  H_2(X)&=\sum_ng\left(\frac{n}{X}\right)\frac{|S(1,1;n)|}{\sqrt{n}}\mu^2(n)\bigg(\frac{k}{2}\bigg)^{\omega(n)}\bigg(\sum_{d|n}\lambda_d\bigg)^2,\\
  H_3(X)&=\sum_ng\left(\frac{n}{X}\right)\frac{S(1,1;n)}{\sqrt{n}}\mu^2(n)\bigg(\sum_{d|n}\lambda_d\bigg)^2.\end{align*}
We wish to estimate
as accurately as possible 
$H_j(X),j=1,2,3$. We shall follow the arguments in \cite{SF2}
and \cite{Ma} to obtain a lower bound for $H_1(X)$. The tools
involved include the Sato-Tate distribution of Kloosterman sums
in prime variables. The investigation on the upper bound for
$H_2(X)$ can be reduced to a problem of evaluating a
multiple-integral, where the Cauchy residue theorem can be
applied. The estimate for $H_3(X)$ is derived
using the spectral
theory of automorphic forms, following
E. Fouvry and Ph. Michel \cite{FM2}.

\begin{proposition}\label{prop:1} For any sufficiently large $X$, we have 
\begin{align*}H_1(X)&\geqslant \widetilde{g}(1)X\mathcal{L}^{-1}(1+o(1))\sum_{2\leqslant i\leqslant5}2^i A_i(\gamma,k)C_i,\end{align*}
where $\gamma$ is defined by $\sqrt{D}=X^\gamma\exp(-\sqrt{\mathcal{L}})$, $A_i(\gamma,k)$ is given by $(\ref{eq:11}),(\ref{eq:12}),(\ref{eq:13}),$ and the constants $C_i$ satisfy $C_2\geqslant 0.11109,C_3\geqslant 0.03557,C_4\geqslant 0.01184,C_5\geqslant 0.00396.$
\end{proposition}

\begin{proposition}\label{prop:2}For any sufficently large $X$, we have 
\begin{align*}H_{2}(X)\leqslant k!^2\cdot R_k(\gamma)\cdot \widetilde{g}(1)X\mathcal{L}^{-1}(1+o(1)),\end{align*}
where $\gamma$ is defined as above and $R$ is a two-variable polynomial given by $(\ref{eq:20}).$\end{proposition}

\begin{proposition}\label{prop:3}For any sufficently large $X$ and $D=O(X^\frac{1}{2}\exp(-\sqrt{\mathcal{L}}))$, we have 
\begin{align*}H_{3}(X)\ll X\mathcal{L}^{-A}\end{align*}for any $A>0.$
\end{proposition}

In order
to obtain a positive lower bound for $H^\pm(X),$
it suffices
to choose $\rho$
so that
\begin{align*}
  \rho H_1(X) > 2 H_2(X) +
  |\rho H_3(X)|
\end{align*}
for $X$ large enough.
For this,
it suffices by the above propositions to choose $k,\gamma$ and $\rho$ so that
\begin{align*}\rho \cdot \sum_{2\leqslant i\leqslant5}2^i A_i(\gamma,k)c_i>2k!^2\cdot R_k(\gamma),\ \ 0<\gamma\leqslant\frac{1}{4}.\end{align*}
With the help of Mathematica 9, we check that the choice
\begin{align*}k=6,\ \ \gamma=\frac{1}{4},\ \ \rho=1.5\times10^5\end{align*}
satisfies the above condition.  We can obtain Theorem \ref{thm:1} since
\begin{align*}\frac{\log\rho}{\log(k/2)}\approx10.849.\end{align*}

\bigskip

\section{Kloosterman sums: From algebraic to analytic}
Kloosterman sums are special kinds of algebraic exponential sums,
which are construced through
algebraic
geometry. Furthermore, Kloosterman sums also appear in the spectral theory of automorphic forms. 
We shall employ both aspects of Kloosterman sums to prove Theorem \ref{thm:1}.

\subsection{Equidistribution of Kloosterman sums: after Katz and Michel}

By the works of Deligne \cite{De} and Katz \cite{Ka2}, the function
$$m\mapsto\frac{S(m,1;p)}{\sqrt{p}}=2\cos\theta_p(m),\ \ m\in \mathbb{F}_p^\times$$
is the Frobenius trace function (restricted to $\mathbf{G}_{m}(\mathbb{F}_p)=\mathbb{F}_p^\times$) of an $\ell$-adic sheaf $\mathcal{K}l$ of rank 2, pure of weight 0 and determinant 1. This means
\begin{align*}2\cos\theta_p(m)=\tr(\Frob_m,\mathcal{K}l).\end{align*}
By the Weyl equidistribution criterion and the Peter-Weyl
theorem, the proof of Katz's equidistribution
theorem
reduces to the study of
\begin{align*}\sum_{m\in\mathbb{F}_p^\times}\sym_k(\theta_p(m))=\sum_{m\in\mathbb{F}_p^\times}\tr(\Frob_m,\sym^k\mathcal{K}l),\end{align*}
where $\sym^k\mathcal{K}l$ is the $k$-th symmetric power
of the Kloosterman sheaf $\mathcal{K}l$ (i.e., the composition of the sheaf $\mathcal{K}l$ with the $k$-th symmetric power representation of $SL_2$) and
\begin{align*}\sym_k(\theta)=\frac{\sin(k+1)\theta}{\sin\theta}.\end{align*}
Using Deligne's main theorem, Katz proved that
\begin{align}\label{eq:6}\left|\sum_{m\in\mathbb{F}_p^\times}\sym_k(\theta_p(m))\right|\leqslant\frac{1}{2}(k+1)p^\frac{1}{2};\end{align}
we refer to
Example 13.6 and the preceding theorem in \cite{Ka2} for more
details.
This implies
that $\left\{ \theta_p(m) : m \in \mathbb{F}_p^\times
\right\}$
equidistributes with respect to the
Sato-Tate
measure $\mu_{\rm ST}$ as $p \rightarrow \infty$.

We can regard (\ref{eq:6}) as the {\it square-root cancellation}
phenomenon for angles of Kloosterman sums.
Due to the supposed randomness of Kloosterman sums,
it is reasonable to expect a similar phenomenon also for $\theta_p(\beta(m))$, where $\beta$ is a non-constant rational function defined over $\mathbb{F}_p^\times$ of fixed degree. 
In fact, we will use this for the map
$\beta:m\mapsto\overline{m}^2$.
In that direction, it is known that
\begin{align}\label{eq:7}\sum_{m\in\mathbb{F}_p^\times}\sym_k(\theta_p(\overline{m}^2))\ll p^\frac{1}{2},\end{align}where the implied constant depends on $k$ polynomially. This estimate has been obtained implicitly by Ph. Michel \cite{Mi}, for whom the relevant sheaf is
\begin{align*}\sym^k([-2]^*\mathcal{K}l).\end{align*}Hence we can conclude from (\ref{eq:7}) the equidistribution of $\theta_p(\overline{m}^2)$ with respect to the Sato-Tate measure $\mu_{\rm ST}$.

\subsection{Equidistribution of Kloosterman sums with composite moduli}
We have been concerned with the equidistribution of Kloosterman sums of prime moduli in the preceding arguments. In later applications, we shall also consider the relevant equidistribution of Kloosterman sums with composite moduli, particularly the products of distinct primes.

Before stating the equidistribution precisely, we would like to
introduce
some measures $\mu^{(j)}$ on $[-1,1]$ which are connected with
the classical Sato-Tate measure $\mu_{\rm ST}$. They
can be defined recursively as follows:
\begin{align*}\ud\mu^{(1)}x=\frac{2}{\pi}\sqrt{1-x^2}\ud x\end{align*}
and \begin{align*}\mu^{(j)}=\mu^{(1)}\otimes\mu^{(j-1)},\ \ \ j\geqslant2.\end{align*}
Then
\begin{align*}\mu^{(1)}([-x,x])=\frac{4}{\pi}\int_0^x\sqrt{1-t^2}\ud t=\frac{2}{\pi}(x\sqrt{1-x^2}+\arcsin x)\end{align*}
and
\begin{align*}\mu^{(j)}([-x,x])=\mu^{(1)}([-x,x])+\frac{4}{\pi}\int_x^1\mu^{(j-1)}([-x/t,x/t])\sqrt{1-t^2}\ud t.\end{align*}

Suppose $p_1,p_2,\cdots,p_k$ are distinct primes. As a
generalization of Katz's result,
one expects that the Kloosterman sum $S(m,1;p_1p_2\cdots p_k)$ equidistributes with respect to the measure $\mu^{(k)}$ as $m$ runs over the primitive residue system and the product $p_1p_2\cdots p_k$ tends to infinity. We shall show this is the case, even while $m$ is restricted to prime variables and the length of sum is sufficiently large compared to the moduli.

\subsection{Equidistribution of Kloosterman sums over prime variables}
We have discussed the equidistribution of Kloosterman sums as the variable runs over consecutive integers in the sense of modular arithmetic. 
Now we consider the equidistribution results when the variable runs amongst the primes, which is in fact what we shall need in the applications to the problem on sign changes. Here we only state the necessary lemmas for the equidistribution, and the precise result will be stated explicitly in the next section.

In a recent series of papers, E. Fouvry, E. Kowalski and
Ph. Michel have investigated analytic properties of some general
functions, known as  {\it algebraic trace functions}, defined over $\mathbb{F}_p$.
In particular, they \cite{FKM} considered the behaviors of such functions over prime variables, and provided a power-saving cancellation. In our applications, we will use the special case of their results.

\begin{lemma}\label{lm:1}
  Let $k$ be a positive integer.
  For each $\varepsilon > 0$
  there exists $\delta=\delta(\varepsilon)>0$
  so that
  if $N>p^{\frac{3}{4}+\varepsilon}$, then
\begin{align*}\sum_{\substack{N<n\leqslant2N\\n\text{ prime}}}\sym_k(\theta_p(\overline{n}^2))\ll Np^{-\delta},\end{align*}
where the implied constant depends on $\varepsilon$ and polynomially on $k$.
\end{lemma}

\proof This is a special case of Theorem 1.5 in \cite{FKM}, where we can take their trace function $K$ as
\begin{align*}n\mapsto\sym_k(\theta_p(\overline{n}^2))\end{align*}defined over $\mathbb{F}_p^\times.$
\endproof

Furthermore, we also require some equidistribution with more than one variables, with respect to the prime moduli and almost prime moduli (products of distinct primes). For the former case, we appeal to the following bilinear form estimate, which can be found in \cite{Mi}, Corollaire 2.11.
\begin{lemma}\label{lm:2}Suppose $1\leqslant M,N\leqslant p.$ For each positive integer $k$ and any coefficients $\alpha=(\alpha_m),\beta=(\beta_n),$ we have
\begin{align*}\mathop{\sum_{M<m\leqslant2M}\sum_{N<n\leqslant2N}}_{(mn,p)=1}\alpha_m\beta_n\sym_k(\theta_p(\overline{mn}^2))\ll 
\|\alpha\|\|\beta\|(MN)^{\frac{1}{2}}(N^{-\frac{1}{2}}+M^{-\frac{1}{2}}p^{\frac{1}{4}}(\log p)^{\frac{1}{2}}),\end{align*}
where $\|\cdot\|$ denotes the $\ell_2$-norm and  the implied constant depends polynomially on $k$.
\end{lemma}

\noindent {\bf Remark.} Lemma \ref{lm:2} is sufficient in our
applications. In fact, we
can remove the restrictions on the sizes of $M,N$ provided that
we insert an extra term
$\|\alpha\|\|\beta\|(MN)^{\frac{1}{2}}p^{-\frac{1}{4}}$ in the
upper bound;
an explicit and more general statement can be found in Theorem 1.17 of \cite{FKM}.

In the case of composite moduli, we require the following estimate, which is stated as Proposition 7.2 in \cite{FM2} and proved by the techniques of $\ell$-adic cohomology.
\begin{lemma}\label{lm:3}Suppose $p_1,p_2,\cdots,p_s$ are
  distinct primes.
  Write $r=p_1p_2\cdots p_s.$ For each $s$-tuple of positive integers $(k_1,k_2,\cdots,k_s),$ and any coefficients $\alpha=(\alpha_m),\beta=(\beta_n),\gamma=(\gamma_{m,n})$ with $m\equiv m'\bmod n\Rightarrow\gamma_{m,n}=\gamma_{m',n},$ we have
\begin{align*}\mathop{\sum_{M<m\leqslant2M}\sum_{N<n\leqslant2N}}_{(mn,r)=1}\alpha_m\beta_n\gamma_{m,n}&\prod_{1\leqslant j\leqslant s}\sym_{k_j}(\theta_{p_j}(\overline{mnrp_j^{-1}}^2))\\
&\ll c(s;\mathbf{k})
\|\alpha\|\|\beta\|\|\gamma\|_{\infty}(MN)^{\frac{1}{2}}(r^{-\frac{1}{8}}+N^{-\frac{1}{4}}r^{\frac{1}{8}}+M^{-\frac{1}{2}}N^{\frac{1}{2}}),\end{align*}
where  $\|\cdot\|_{\infty}$ denotes the sup-norm, $c(s;\mathbf{k})=3^s\prod_{j=1}^s(k_j+1)$ and the implied constant is absolute.
\end{lemma}

\bigskip

\section{Proof of Proposition \ref{prop:1}: Lower bound for $H_1(X)$}

\subsection{Initial step: Preparation for equidistribution}
We start the proof of Proposition \ref{prop:1}. 
Let
\[C(m,n)=\frac{S(\overline{m}^2,1;n)}{2^{\omega(n)}\sqrt{n}}\]for
$(m,n)=1$. Then we have $|C(m,n)|\leqslant1$ for squarefree $n$
by (\ref{eq:2}), and it follows from
the
Chinese remainder theorem that
\begin{align}\label{eq:8}C(1,mn)=C(m,n)C(n,m).\end{align}
In particular, we have $C(m,p)=\cos\theta_p(\overline{m}^2)$ for $(m,p)=1$. In this way, we have
\begin{align*}H_1(X)&= \sum_{n}g\left(\frac{n}{X}\right)\mu^2(n)2^{\omega(n)}|C(1,n)|\bigg(\sum_{d|n}\lambda_d\bigg)^2.\end{align*}
In our applications,
we
need only consider
those $n$ with few prime factors. To that end, we introduce the interval
\begin{align*}I(P)=(P,P+P\mathcal{L}^{-1}],\end{align*}
and the set of the products of primes
\begin{align*}\mathcal{P}_i(X;P_{i1},P_{i2},\cdots,P_{ii})&=\{p_1p_2\cdots p_i:p_j\in I(P_{ij})\text{ for each } j\leqslant i\}\end{align*}
for each positive integer $i\geqslant2.$ Furthermore, for each
fixed $i$, we assume that $\{P_{ij}\}$ is a decreasing sequence
as $j$ varies and the product of the lengths
of the intervals $I(P_{ij})$ is exactly $X$, i.e.,
that
\begin{align}\label{eq:9}P_{i1}>P_{i2}>\cdots>P_{ii}>X^\varepsilon,\ \ \prod_{1\leqslant j\leqslant i}|I(P_{ij})|=X.\end{align}
In this way, we can bound $H_1(X)$ from below by the summation
over $\mathcal{P}_i(X;P_{i1},P_{i2},\cdots,P_{ii})$; for this, we employ the variants of 
the Sato-Tate distributions stated above.
Due to the positivity of each term, we can drop those $n$'s with
``bad''
arithmetic structures. To this end, we introduce the
following
restrictions on the size of $P_{ij}$:
\begin{equation}\label{eq:10}\begin{split}\begin{cases}P_{21}^{3/4}X^\eta<P_{22},\ \ \eta=10^{-2014},\\
P_{31}^{1/2}\exp(\sqrt{\mathcal{L}})<P_{32},\\
P_{41}^{1/2}\exp(\sqrt{\mathcal{L}})<P_{42}P_{43},\\
P_{51}^{1/2}\exp(\sqrt{\mathcal{L}})<P_{52}P_{53}P_{54}\text{ and } (P_{53}P_{54}P_{55})^{1/2}\exp(\sqrt{\mathcal{L}})<P_{52},\\
\cdots\end{cases}\end{split}\end{equation}
Now summing up to $i=5,$ we have the lower bound
\begin{align*}H_1(X)&\geqslant\sum_{2\leqslant i\leqslant5}2^iH_{1,i}(X),\end{align*}
where
\begin{align*}H_{1,i}(X)&=\sideset{}{^\dagger}\sum_{P_{i1},P_{i2},\cdots,P_{ii}}\ \ \sum_{n\in\mathcal{P}_i(X;P_{i1},P_{i2},\cdots,P_{ii})}g\left(\frac{n}{X}\right)|C(1,n)|\bigg(\sum_{d|n}\lambda_d\bigg)^2\end{align*}
with
the symbol $\dagger$ denoting the restrictions (\ref{eq:9}) and (\ref{eq:10}).

Recalling the choice (\ref{eq:4}), we find, for each $n\in\mathcal{P}_i(X;P_{i1},P_{i2},\cdots,P_{ii})$, that
\begin{align*}\sum_{d|n}\lambda_d=(1+o(1))L_i(\gamma,k;X^{\alpha_1},X^{\alpha_2},\cdots,X^{\alpha_i}),\end{align*}
where
\begin{align}\label{eq:11}L_i(\gamma,k;X^{\alpha_1},X^{\alpha_2},\cdots,X^{\alpha_i})=\sum_{\substack{\mathcal{A}\subseteq\{\alpha_1,\alpha_2,\cdots,\alpha_i\}\\\sum_{\alpha\in\mathcal{A}}\alpha<\gamma}}(-1)^{|\mathcal{A}|}\bigg(1-\frac{1}{\gamma}\sum_{\alpha\in\mathcal{A}}\alpha\bigg)^k.\end{align}

Note the bound $|C(1,n)|\leqslant1$.
From partial summation, we can write
\begin{align*}H_{1,i}(X)&=\widetilde{g}(1)\mathcal{L}^{i-1}(1+o(1))\int\cdots\int_{\mathcal{R}_i} L_i^2(\gamma,k;X^{1-\alpha_2-\cdots-\alpha_i},X^{\alpha_2},\cdots,X^{\alpha_i})\ud\alpha_2\cdots\ud\alpha_i\\
&\ \ \ \ \ \times\sum_{n\in\mathcal{P}_i(X;X^{1-\alpha_2-\cdots-\alpha_i},X^{\alpha_2},\cdots,X^{\alpha_i})}|C(1,n)|,\end{align*}
where the multiple-integral is over the area $\mathcal{R}_i$:
\begin{equation}\label{eq:12}\begin{split}
\mathcal{R}_2&:=\{\alpha_2\in(0,1):(\frac{3}{4}+\eta)(1-\alpha_2)<\alpha_2<\frac{1}{2}\},\ \ \eta=10^{-2014},\\
\mathcal{R}_3&:=\{(\alpha_2,\alpha_3)\in(0,1)^2:\frac{1}{2}(1-\alpha_2-\alpha_3)<\alpha_2,\alpha_3<\alpha_2<1-\alpha_2-\alpha_3\},\\
\mathcal{R}_4&:=\{(\alpha_2,\alpha_3,\alpha_4)\in(0,1)^3:\frac{1}{2}(1-\alpha_2-\alpha_3-\alpha_4)<\alpha_2+\alpha_3\}\\
&\ \ \ \ \ \cap\{(\alpha_2,\alpha_3,\alpha_4)\in(0,1)^3:\alpha_4<\alpha_3<\alpha_2<1-\alpha_2-\alpha_3-\alpha_4\},\\
\mathcal{R}_5&:=\{(\alpha_2,\alpha_3,\alpha_4,\alpha_5)\in(0,1)^4:\frac{1}{2}(1-\alpha_2-\alpha_3-\alpha_4-\alpha_5)<\alpha_2+\alpha_3+\alpha_4\}\\
&\ \ \ \ \ \cap\{(\alpha_2,\alpha_3,\alpha_4,\alpha_5)\in(0,1)^4:\frac{1}{2}(\alpha_3+\alpha_4+\alpha_5)<\alpha_2\}\\
&\ \ \ \ \ \cap\{(\alpha_2,\alpha_3,\alpha_4,\alpha_5)\in(0,1)^4:\alpha_5<\alpha_4<\alpha_3<\alpha_2<1-\alpha_2-\alpha_3-\alpha_4-\alpha_5\}.\end{split}\end{equation}

\subsection{Applications of equidistribution}

In the preceding ranges, we deduce from Lemmas \ref{lm:1},
\ref{lm:2} and \ref{lm:3} the following equidistribution
results, which extend Propositions 6.1, 6.2 and 6.3 in \cite{FM2}.
\begin{lemma}\label{lm:4}
  With the notation as above, for $i\in\{2,3,4,5\}$ and $(\alpha_2,\cdots,\alpha_i)\in\mathcal{R}_i$, the sets
  \[\{C(p_1,p_2\cdots p_j):n=p_1p_2\cdots p_j
  \in\mathcal{P}_i(X;X^{1-\alpha_2-\cdots-\alpha_i},X^{\alpha_2},\cdots,X^{\alpha_i})\}\]
and
\[\{C(p_2\cdots p_j,p_1):n=p_1p_2\cdots p_j
\in\mathcal{P}_i(X;X^{1-\alpha_2-\cdots-\alpha_i},X^{\alpha_2},\cdots,X^{\alpha_i})\}\]
equidistribute in $[-1,1]$ with respect to $\mu^{(i-1)}$ and
$\mu^{(1)}$,
respectively, as $X\rightarrow+\infty$.
\end{lemma}

Lemma \ref{lm:4} provides the equidistribution
of Kloosterman sums with fixed moduli. However, for the purpose
of lower bound for $H_{1,i}(X)$, we must understand the
distribution of $C(1,n)$ as $n$ runs over
$\mathcal{P}_i(X;X^{1-\alpha_2-\cdots-\alpha_i},X^{\alpha_2},\cdots,X^{\alpha_i})$. Of
course, this would partially follow from the factorization of
$C(1,n)$ as the product the two Kloosterman sums, of which we
know equidistribution in the ranges stated in Lemma
\ref{lm:4}. Hence, in general, we are
faced with the problem of obtaining the result for joint distribution of
two sequences assuming equidistribution of each. For this, we appeal to the following rearrangement type inequality due to K. Matom\"{a}ki \cite{Ma}.
\begin{lemma}\label{lm:5}Assume that the sequences $(a_n)_{n\leqslant N}$ and $(b_n)_{n\leqslant N}$ contained
  in $[0,1]$  equidistribute with respect to some
  absolutely continuous measures $\mu_a$ and $\mu_b$, respectively, as $N\rightarrow\infty$. Then
\begin{align*}(1+o(1))\int_0^1xy_l(x)\ud\mu_a([0,x]) \leqslant\frac{1}{N}\sum_{n\leqslant N}a_nb_n\leqslant
(1+o(1))\int_0^1xy_u(x)\ud\mu_a([0,x]),\end{align*}
where $y_l(x)$ is the smallest solution to the equation $\mu_b([y_l,1])=\mu_a([0,x])$ and
$y_u(x)$ is the largest solution to the equation $\mu_b([0,y_u])=\mu_a([0,x])$.
\end{lemma}

Now we write
\begin{align*}&\sum_{n\in\mathcal{P}_i(X;X^{1-\alpha_2-\cdots-\alpha_i},X^{\alpha_2},\cdots,X^{\alpha_i})}|C(1,n)|\\
&=\sum_{n\in\mathcal{P}_i(X;X^{1-\alpha_2-\cdots-\alpha_i},X^{\alpha_2},\cdots,X^{\alpha_i})}|C(p_2\cdots p_i,p_1)||C(p_1,p_2\cdots p_i)|.\end{align*}

By Lemmas \ref{lm:4} and \ref{lm:5}, this is
\begin{align*}\geqslant |\mathcal{P}_i(X;X^{1-\alpha_2-\cdots-\alpha_i},X^{\alpha_2},\cdots,X^{\alpha_i})|C_i(1+o(1))=\frac{X\mathcal{L}^{-i}C_i(1+o(1))}{\alpha_2\cdots\alpha_j(1-\alpha_2-\cdots-\alpha_j)}\end{align*}
for some positive constant $C_i.$ Hence we can obtain the inequality
\begin{align*}H_1(X)&\geqslant \widetilde{g}(1)X\mathcal{L}^{-1}(1+o(1))\sum_{2\leqslant i\leqslant5}2^i A_i(\gamma,k)C_i,\end{align*}
where
\begin{align}\label{eq:13}A_i(\gamma,k)&=\int\cdots\int_{\mathcal{R}_i} \frac{L_i^2(\gamma,k;X^{1-\alpha_2-\cdots-\alpha_i},X^{\alpha_2},\cdots,X^{\alpha_i})}{\alpha_2\cdots\alpha_j(1-\alpha_2-\cdots-\alpha_j)}\ud\alpha_2\cdots\ud\alpha_i.\end{align}
More precisely, by Lemma \ref{lm:5}, we can take
\begin{align*}C_i\geqslant\int_0^1x y_i(x)\ud\mu^{(1)}([-x,x]),\end{align*}
where $y_i(x)$ is the unique solution to the equation
\begin{align*}\mu^{(1)}([-x,x])=\mu^{(i-1)}([-1,-y]\cup[y,1])=1-\mu^{(i-1)}([-y,y]).\end{align*}
With the help of Mathematica 9, we can obtain
\begin{align*}
C_2&\geqslant 0.11109,\\
C_3&\geqslant 0.03557,\\
C_4&\geqslant 0.01184,\\
C_5&\geqslant 0.00396.
\end{align*}
This proves Proposition \ref{prop:1}.

\bigskip

\section{Proof of Proposition \ref{prop:2}: Upper bound for $H_2(X)$}
Before starting the proof of Proposition \ref{prop:2}, we state two results from complex analysis. The first one is an example of the Mellin inversion formula, as an immediate consequence of Cauchy's residue theorem.

\begin{lemma}\label{lm:6}Suppose $k$ is a non-negative integer. For any positive number $x$, we have
\begin{align*}\frac{k!}{2\pi i}\int_{1-i\infty}^{1+i\infty}\frac{x^s}{s^{k+1}}\ud s=\begin{cases}0,\ \ \ &0<x\leqslant1,\\
(\log x)^k,&x>1.\end{cases}\end{align*}
\end{lemma}

The following lemma is contained implicitly in \cite{GMPY}.
\begin{lemma}\label{lm:7} Suppose $x\geqslant1,$ and $k,l$ are non-negative integers. Then we have
\begin{align*}\Res_{(s_1,s_2)=(0,0)}\frac{x^{s_1+s_2}}{(s_1+s_2)^l(s_1s_2)^{k+1}}=\frac{1}{(2k+l)!}\binom{2k}{k}(\log x)^{2k+l}.\end{align*}
\end{lemma}
\proof We adopt the method of Y. Motohashi,
as presented in \cite{GMPY}. Consider the double-integral
\begin{align*}J=\frac{1}{(2\pi i)^2}\int_{\mathcal{D}_1}\int_{\mathcal{D}_2}\frac{x^{s_1+s_2}}{(s_1+s_2)^l(s_1s_2)^{k+1}}\ud s_1\ud s_2,\end{align*}
where $\mathcal{D}_1,\mathcal{D}_2$ are two small circles
centered at
the origin of radii $\varepsilon,2\varepsilon$, respectively. One can see that $J$ just represents the residue of the integrands at origin. Write $s_1=s,s_2=s\xi,$ then $J$ becomes
\begin{align*}J=\frac{1}{(2\pi i)^2}\int_{\mathcal{D}_1}\int_{\mathcal{D}_3}\frac{x^{(\xi+1)s}}{s^{2k+l+1}\xi^{k+1}(\xi+1)^l}\ud s\ud \xi,\end{align*}
where $\mathcal{D}_3$ is the circle centered at origin of radius $2.$ Clearly, the $s$-integral is
\begin{align*}\frac{1}{(2k+l)!}((\xi+1)\log x)^{2k+l}.\end{align*}
Now it follows that
\begin{align*}J=\frac{(\log x)^{2k+l}}{(2k+l)!}\frac{1}{2\pi i}\int_{\mathcal{D}_3}\frac{(\xi+1)^{2k}}{\xi^{k+1}}\ud \xi=\frac{(\log x)^{2k+l}}{(2k+l)!}\binom{2k}{k}\end{align*} since the $\xi$-integral detects the coefficient of $\xi^k$ in the expansion of $(1+\xi)^{2k}.$
This establishes Lemma  \ref{lm:7}.
\endproof

\subsection{Expressing as a multiple-integral}
The arguments in this section have almost nothing to do with
Kloosterman sums;
the only fact we shall use is that $|C(1,n)|\leqslant1.$
More precisely, we have
\begin{align*}H_2(X)&=\sum_{n}g\left(\frac{n}{X}\right)|C(1,n)|\mu^2(n)k^{\omega(n)}\bigg(\sum_{d|n}\lambda_d\bigg)^2\\
&\leqslant\sum_{n}g\left(\frac{n}{X}\right)\mu^2(n)k^{\omega(n)}\sum_{d|n}\xi(d),\end{align*}
where
\begin{align}\label{eq:14}\xi(d)=\sum_{[d_1,d_2]=d}\lambda_{d_1}\lambda_{d_2}.\end{align}
Furthermore, we have
\begin{align*}H_2(X)&\leqslant\sum_{d}\xi(d)k^{\omega(d)}\sum_{(n,d)=1}g\left(\frac{nd}{X}\right)\mu^2(n)k^{\omega(n)}.\end{align*}
Now we would like to evaluate the $n$-sum. By Mellin inversion, we can write
\begin{align*}g\left(\frac{nd}{X}\right)=\frac{1}{2\pi i}\int_{2-i\infty}^{2+i\infty}\widetilde{g}(s)\left(\frac{X}{nd}\right)^s\ud s.\end{align*}Then it follows that
\begin{align*}\sum_n=\frac{1}{2\pi i}\int_{2-i\infty}^{2+i\infty}\widetilde{g}(s)\left(\frac{X}{d}\right)^sT(d,s)\ud s,\end{align*}
where $T(d,s)$ is defined by the Dirichlet series
\begin{align*}T(d,s)=\sum_{\substack{n\geqslant1\\(n,d)=1}}\frac{\mu^2(n)k^{\omega(n)}}{n^s},\ \ \Re s>1.\end{align*}
For $\Re s>1,$ we have
\begin{align*}T(d,s)=\prod_{p\nmid d}\bigg(1+\frac{k}{p^s}\bigg)=\zeta^k(s)T^*(d,s),\end{align*}
where, for each fixed positive integer $d$, $T^*(d,s)$ is a holomorphic function in the half plane $\Re s>0$ and
$$T^*(d,1)=\prod_{p\mid d}\bigg(1+\frac{k}{p}\bigg)^{-1}\cdot\prod_p\bigg(1+\frac{k}{p}\bigg)\bigg(1-\frac{1}{p}\bigg)^k.$$ Thus $T(d,s)$ admits a meromorphic continuation to $\Re s>0$ with $s=1$ as the unique pole of order $k$. Note that
\begin{align*}\sum_n=\frac{1}{2\pi i}\int_{2-i\infty}^{2+i\infty}\widetilde{g}(s)\left(\frac{X}{d}\right)^s\zeta^k(s)T^*(d,s)\ud s.\end{align*}Moving the integral line to $\Re s=\frac{1}{2}$, we shall pass the pole $s=1$, getting
\begin{align*}\sum_n=\Res_{s=1}\widetilde{g}(s)\left(\frac{X}{d}\right)^s\zeta^k(s)T^*(d,s)+\frac{1}{2\pi i}\int_{\frac{1}{2}-i\infty}^{\frac{1}{2}+i\infty}\widetilde{g}(s)\left(\frac{X}{d}\right)^s\zeta^k(s)T^*(d,s)\ud s.\end{align*}
From the growth of the integrand, we can easily verify that the
second term
is
bounded by $(X/d)^\delta$ for some $\delta<1.$ The first term is in fact
\begin{align*}\frac{\widetilde{g}(1)X}{d}&T^*(d,1)\Res_{s=1}\zeta^k(s)\\&=\frac{\widetilde{g}(1)}{(k-1)!}\prod_{p}\bigg(1-\frac{1}{p}\bigg)^{k}\bigg(1+\frac{k}{p}\bigg)\cdot\prod_{p|d}\bigg(1+\frac{k}{p}\bigg)^{-1}\frac{X}{d}P_{k-1}(\log(X/d)),\end{align*}
where $P_k(\cdot)$ is a monic polynomial of degree $k-1$.

Hence we find
\begin{align*}\sum_n=\frac{\widetilde{g}(1)}{(k-1)!}\prod_{p}\bigg(1-\frac{1}{p}\bigg)^{k}\bigg(1+\frac{k}{p}\bigg)\cdot\prod_{p|d}\bigg(1+\frac{k}{p}\bigg)^{-1}\frac{X}{d}P_{k-1}(\log(X/d))+O(Xd^{-1}\mathcal{L}^{-A}),\end{align*}
and it follows that
\begin{align}\label{eq:15}H_2(X)&\leqslant\frac{\widetilde{g}(1)}{(k-1)!}\prod_{p}\bigg(1-\frac{1}{p}\bigg)^{k}\bigg(1+\frac{k}{p}\bigg)X\cdot D(X)+O(X\mathcal{L}^{-A}),\end{align}
where
\begin{align}\label{eq:16}D(X)&=\sum_{d\leqslant X}\xi(d)\frac{k^{\omega(d)}}{d}\prod_{p|d}\bigg(1+\frac{k}{p}\bigg)^{-1}P_{k-1}(\log(X/d)).\end{align}

Denote by $\widetilde{D}(X)$ the relevant contribution from the highest order monomial in $P_{k-1}(\log(X/d))$.
Then by Lemma \ref{lm:6} we have
\begin{align}\label{eq:17}\widetilde{D}(X)&=\frac{(k-1)!}{2\pi i}\int_{1-i\infty}^{1+i\infty}M(s)\frac{X^s}{s^k}\ud s,\end{align}
where
\begin{align*}M(s)&=\sum_{d\geqslant 1}\xi(d)\frac{k^{\omega(d)}}{d^{s+1}}\prod_{p|d}\bigg(1+\frac{k}{p}\bigg)^{-1}.\end{align*}

Now define
\begin{align*}u(n,s)=\frac{k^{\omega(n)}}{n^{s+1}}\prod_{p|n}\bigg(1+\frac{k}{p}\bigg)^{-1},\end{align*}
and rewrite $M(s)$ as
\begin{align}\label{eq:18}M(s)&=\sum_{d\geqslant 1}\xi(d)u(d,s)=\mathop{\sum\sum}_{d_1,d_2\leqslant \sqrt{D}}\lambda_{d_1}\lambda_{d_2}u([d_1,d_2],s).\end{align}
Note that 
\begin{align*}u([d_1,d_2],s)&=u(d_1,s)u(d_2,s)\bigg(\prod_{p|(d_1,d_2)}\frac{1}{u(p,s)}\bigg)\\
&=u(d_1,s)u(d_2,s)\sum_{d|(d_1,d_2)}\prod_{p|d}\bigg(\frac{1}{u(p,s)}-1\bigg),\end{align*}
thus (\ref{eq:18}) becomes
\begin{align}\label{eq:19}M(s)=\sum_{m\leqslant \sqrt{D}}\mu^2(m)\prod_{p|m}\bigg(\frac{1}{u(p,s)}-1\bigg)Z(m,s)^2,\end{align}
where
\begin{align*}Z(m,s)&=\frac{1}{(\log \sqrt{D})^k}\sum_{\substack{d\leqslant \sqrt{D}\\m|d}}\mu(d)u(d,s)\log^k(\sqrt{D}/d)\\
&=\frac{1}{(\log\sqrt{D})^k}\mu(m)u(m,s)\sum_{\substack{d\leqslant \sqrt{D}/m\\(d,m)=1}}\mu(d)u(d,s)\log^k(\sqrt{D}/md).\end{align*}

Write
\begin{align*}G(w,m)&=\sum_{\substack{d\geqslant1\\(d,m)=1}}\mu(d)u(d,w),\end{align*}
for $\Re w>0.$ It is clear that $G(w,m)$ has an analytic continuation on $\mathbb{C}$ in the $w$-variable, in fact, we can write
\begin{align*}G(w,m)&=\frac{1}{\zeta(w+1)^k}F(w)\cdot\prod_{p|m}(1-u(p,w))^{-1},\end{align*}
where $F(w)$ is defined by
\begin{align*}F(w)&=\prod_{p}\bigg(1-\frac{k}{p^{w}(p+k)}\bigg)\bigg(1-\frac{1}{p^{w+1}}\bigg)^{-k}\end{align*}as $\Re w>-1.$ Now we have
\begin{align*}Z(m,s)&=\frac{\mu(m)u(m,s)}{(\log\sqrt{D})^k}\frac{k!}{2\pi i}\int_{1-i\infty}^{1+i\infty}G(w+s;m)\frac{(\sqrt{D}/m)^w}{w^{k+1}}\ud w,\end{align*}
from which and (\ref{eq:19}) we find
\begin{align*}M(s)&=\frac{k!^2}{(\log\sqrt{D})^{2k}}\frac{1}{(2\pi i)^2}\int_{1-i\infty}^{1+i\infty}\int_{1-i\infty}^{1+i\infty}\frac{F(w_1+s)F(w_2+s)}{\zeta^k(w_1+s+1)\zeta^k(w_2+s+1)}\frac{(\sqrt{D})^{w_1+w_2}}{(w_1w_2)^{k+1}}\ud w_1\ud w_2\\
&\ \ \ \ \times\sum_{m\leqslant \sqrt{D}}\frac{\mu^2(m)u^2(m,s)}{m^{w_1+w_2}}\prod_{p|m}\bigg(\frac{1}{u(p,s)}-1\bigg)(1-u(p,s+w_1))^{-1}(1-u(p,s+w_2))^{-1}.\end{align*}
Note that the $m$-sum can be expressed as
\begin{align*}\frac{1}{2\pi i}\int_{1-i\infty}^{1+i\infty}\zeta^k(t+w_1+w_2+s+1)H(t,w_1,w_2,s)\frac{(\sqrt{D})^t}{t}\ud t,\end{align*}
where $H(t,w_1,w_2,s)$ is holomorphic for $\Re(t+s),\Re(w_1+s),\Re(w_2+s)>-1$ with
\begin{align*}H(0,0,0,0)=\prod_{p}\bigg(1+\frac{k}{p}\bigg)\bigg(1-\frac{1}{p}\bigg)^k.\end{align*}
Hence we can deduce from (\ref{eq:17}) that
\begin{align*}\widetilde{D}(X)&=\frac{(k-1)!\cdot k!^2}{(\log\sqrt{D})^{2k}}\frac{1}{(2\pi i)^4}\iiiint\limits_{(1)(1)(1)(1)}K(t,w_1,w_2,s)\\
&\ \ \ \ \ \ \times\bigg(\frac{(w_1+s)(w_2+s)}{s(t+w_1+w_2+s)}\bigg)^k\frac{(\sqrt{D})^{t+w_1+w_2}X^s}{t(w_1w_2)^{k+1}}\ud t\ud w_1\ud w_2\ud s,\end{align*}
where $K(t,w_1,w_2,s)$ is holomorphic for $\Re t,\Re w_1,\Re w_2,\Re s,\Re(t+s),\Re(w_1+s),\Re(w_2+s)>-1$ and
\begin{align*}K(0,0,0,0)=\prod_{p}\bigg(1+\frac{k}{p}\bigg)^{-1}\bigg(1-\frac{1}{p}\bigg)^{-k}.\end{align*}

\subsection{Shifting contours}
Now we are in the position to evaluate the multiple-integral by shifting contours. To this end, we define 
\begin{align*}\mathcal{C}=\{-\frac{1}{2014\log(|t|+2)}+it:t\in\mathbb{R}\},\end{align*}
which is related to the zero-free region of Riemann zeta functions.

We can shift all contours to $\sigma=1/\log X$ without passing any poles of the integrand. We now continue to shift the four contours to $\mathcal{C}$ one by one; we consider the $t$-integral first. There are two singularities $t=0$ and $t=-(w_1+w_2+s)$, which are of multiplicity 1 and $k$, respectively.
Hence, after the shifting, the new integrand becomes
\begin{align*}&K(0,w_1,w_2,s)\bigg(\frac{(w_1+s)(w_2+s)}{s(w_1+w_2+s)}\bigg)^k\frac{(\sqrt{D})^{w_1+w_2}X^s}{(w_1w_2)^{k+1}}+K(-(w_1+w_2+s),w_1,w_2,s)\\
&\ \ \ \ \ \times\bigg(\frac{(w_1+s)(w_2+s)}{s}\bigg)^k\frac{(\sqrt{D})^{w_1+w_2}X^s}{(w_1w_2)^{k+1}}\frac{1}{(k-1)!}\frac{\partial^{k-1}}{\partial t^{k-1}}\frac{(\sqrt{D})^t}{t}\bigg|_{t=-(w_1+w_2+s)}.\end{align*}
In fact, there is also another contribution from the integal along $\mathcal{C}$, which is of a lower order of magnitute due to the growth of Riemann zeta functions (In the discussion below, we shall not present explicitly the error terms resulting from shifting contours). Note that the second term comes from the singularity $t=-(w_1+w_2+s),$ and the factor $(\sqrt{D})^{w_1+w_2}$ will vanish after taking the partial derivatives, thus we conclude from Lemma \ref{lm:6} that the second term will produce a contribution of lower order of magnitude. We only consider the first term in latter discussions since what we are interested in is the constant in the main term.

Now we are left with the triple-integral with respect to $w_1,w_2$ and $s$. The resulting integrand is
\begin{align*}K(0,w_1,w_2,s)\bigg(\frac{(w_1+s)(w_2+s)}{s(w_1+w_2+s)}\bigg)^k\frac{(\sqrt{D})^{w_1+w_2}X^s}{(w_1w_2)^{k+1}}.\end{align*}
Now we turn to shift the $s$-contour. Clearly, we shall
encontour four singularities $s=0,-w_1,-w_2$ and
$-(w_1+w_2)$. In fact, the latter three ones will produce
factors of the shape $(\sqrt{D}/X)^{w_1},(\sqrt{D}/X)^{w_2}$ and
$(\sqrt{D}/X)^{w_1+w_2}$. Following the same arguments as above,
we conclude from Lemma \ref{lm:6} that all of these will
contribute
negligibly. Hence we need only consider the singularity $s=0.$ Note that
\begin{align*}\bigg(\frac{(w_1+s)(w_2+s)}{s(w_1+w_2+s)}\bigg)^k=\bigg(1+\frac{w_1w_2}{s(w_1+w_2+s)}\bigg)^k=\sum_{j=0}^{k}\binom{k}{j}\bigg(\frac{w_1w_2}{s(w_1+w_2+s)}\bigg)^{j},\end{align*}thus we can rewrite the integrand as
\begin{align*}\sum_{j=0}^{k}\binom{k}{j}K_j(w_1,w_2,s),\end{align*}
where
\begin{align*}K_j(w_1,w_2,s)=K(0,w_1,w_2,s)\bigg(\frac{w_1w_2}{s(w_1+w_2+s)}\bigg)^{j}\frac{(\sqrt{D})^{w_1+w_2}X^s}{(w_1w_2)^{k+1}}.\end{align*}
For $j\geqslant1,$ we have
\begin{align*}\Res_{s=0}K_j(w_1,w_2,s)&=\frac{1}{\Gamma(j)}\frac{(\sqrt{D})^{w_1+w_2}}{(w_1w_2)^{k+1-j}}\bigg(\frac{\partial^{j-1}}{\partial s^{j-1}}\frac{X^s}{(w_1+w_2+s)^j}\bigg)_{s=0}\\
&=\frac{1}{\Gamma(j)}\sum_{i=0}^{j-1}\binom{j-1}{i}(\log X)^{j-i-1}(-1)^i\frac{\Gamma(j+i)}{\Gamma(j)}\\
&\ \ \ \ \times\frac{(\sqrt{D})^{w_1+w_2}}{(w_1w_2)^{k+1-j}(w_1+w_2)^{j+i}}.\end{align*}
Repeating the same arguments to the $w_1,w_2$-integrals, it
follows that we need only
consider the residue at $w_1=w_2=0$, thus we deduce from Lemma \ref{lm:7} that
\begin{align*}\Res_{(0,0,0)}K_j(w_1,w_2,s)&=\sum_{i=0}^{j-1}\binom{j-1}{i}\binom{2(k-j)}{k-j}\frac{\Gamma(j+i)}{\Gamma(j)^2}\frac{(-1)^i}{(2k-j+i)!}\\
&\ \ \ \ \times(\log X)^{j-i-1}(\log \sqrt{D})^{2k-j+i}.\end{align*}

Hence we obtain
\begin{align*}\widetilde{D}(X)&=(k-1)!\cdot k!^2(1+o(1))\prod_{p}\bigg(1+\frac{k}{p}\bigg)^{-1}\bigg(1-\frac{1}{p}\bigg)^{-k}\sum_{j=1}^{k}\sum_{i=0}^{j-1}\binom{k}{j}\binom{j-1}{i}\binom{2(k-j)}{k-j}\\
&\ \ \ \ \times\frac{\Gamma(j+i)}{\Gamma(j)^2}\frac{(-1)^i}{(2k-j+i)!}(\log X)^{j-i-1}(\log \sqrt{D})^{-j+i}.\end{align*}

\subsection{Conclusion}
By similar arguments, we can obtain
an
asymptotic formula for the contributions related to lower order
terms of the shape $P_{k-1}(\log(X/d))$.
Comparing with the above asymptotic formula for
$\widetilde{D}(X)$, we find that $\widetilde{D}(X)$ contributes the
main term in (\ref{eq:16}).
It then follows from (\ref{eq:15}) that
\begin{align*}H_{2}(X)&\leqslant k!^2\cdot\widetilde{g}(1)X\mathcal{L}^{-1}(1+o(1))R_k(\gamma),\end{align*}
where $\gamma$ is defined by $\sqrt{D}=X^\gamma\exp(-\sqrt{\mathcal{L}})$ and
\begin{align}\label{eq:20}R_k(y)=\sum_{j=1}^{k}\sum_{i=0}^{j-1}\binom{k}{j}\binom{j-1}{i}\binom{2(k-j)}{k-j}\frac{\Gamma(j+i)}{\Gamma(j)^2}\frac{(-1)^i}{(2k-j+i)!}\frac{1}{y^{j-i}}.\end{align}
This completes the proof of Proposition \ref{prop:2}.

\bigskip

\section{Proof of Proposition \ref{prop:3}: Estimate for $H_3(X)$}
Opening the square in $H_3(X)$ and switching the summations, we get
\begin{align*}H_3(X)&=\sum_{d\leqslant D}\xi(d)\sum_{n\equiv0\bmod d}\mu^2(n)g\left(\frac{n}{X}\right)\frac{S(1,1;n)}{\sqrt{n}},\end{align*}
where $\xi(d)$ is defined by (\ref{eq:14}), giving $|\xi(d)|\leqslant 3^{\omega(d)}$ for any squarefree $d$. Hence we have
\begin{align*}H_3(X)&\ll\sum_{d\leqslant D}3^{\omega(d)}\left|\sum_{n\equiv0\bmod d}\mu^2(n)g\left(\frac{n}{X}\right)\frac{S(1,1;n)}{\sqrt{n}}\right|.\end{align*}

Now we are in a position to estimate
mean values
of Kloosterman sums.
We appeal to the following Bombieri-Vinogradov type theorem for
Kloosterman sums, which has been proved in \cite{FM2} using the
spectral theory of automorphic forms without the extra factor
$\mu^2(n)$;
the version employed here
is due to \cite{SF2} as Corollaire 2.2 therein.
\begin{lemma}\label{lm:8}For any $A>0$ there exists some $B=B(A)>0$ such that
\begin{align*}\sum_{q\leqslant \sqrt{X}\mathcal{L}^{-B}}\left|\sum_{n\equiv0\bmod q}\mu^2(n)g\left(\frac{n}{X}\right)\frac{S(1,1;n)}{\sqrt{n}}\right|\ll X\mathcal{L}^{-A},\end{align*}where the implied constant depends on $A$ and $g.$
\end{lemma}

Proposition \ref{prop:3} can be established by the following lemma, which is weighted by divisor functions.
\begin{lemma}\label{lm:9}For any $A>0$, there exists some $B=B(A)>0$ such that
\begin{align*}\sum_{q\leqslant \sqrt{X}\mathcal{L}^{-B}}3^{\omega(q)}\left|\sum_{n\equiv0\bmod q}\mu^2(n)g\left(\frac{n}{X}\right)\frac{S(1,1;n)}{\sqrt{n}}\right|\ll X\mathcal{L}^{-A},\end{align*}where the implied constant depends on $A$ and $g.$
\end{lemma}

\proof
For any fixed $A>0$, we split the $q$-sum as
\begin{align*}\sum_{3^{\omega(q)}\leqslant \mathcal{L}^{A/2}}+\sum_{3^{\omega(q)}> \mathcal{L}^{A/2}},\end{align*}
hence the contribution from the first term is at most $O(X\mathcal{L}^{-A/2})$ by Lemma \ref{lm:8}. For the second term, the contribution is
\begin{align*}&\ll\mathcal{L}^{-A/2}\sum_{q\leqslant \sqrt{X}\mathcal{L}^{-B}}\mu^2(q)9^{\omega(q)}\left|\sum_{n\equiv0\bmod q}\mu^2(n)g\left(\frac{n}{X}\right)\frac{S(1,1;n)}{\sqrt{n}}\right|\\
&\ll\mathcal{L}^{-A/2}\sum_{q\leqslant \sqrt{X}\mathcal{L}^{-B}}\mu^2(q)9^{\omega(q)}\sum_{\substack{n\sim X\\n\equiv0\bmod q}}2^{\omega(n)}\\
&\ll X\mathcal{L}^{1-A/2}\sum_{q\leqslant \sqrt{X}\mathcal{L}^{-B}}\mu^2(q)\frac{18^{\omega(q)}}{q}\\
&\ll X\mathcal{L}^{18-A/2}.\end{align*}
Now the lemma follows from the arbitrariness of $A$.\endproof

\bigskip

\bigskip

\bibliographystyle{plainnat}

\begin{thebibliography}{123}
\bibitem{De}P. Deligne. La conjecture de Weil II. \emph{Publ. Math. IHES} \textbf{52} (1980), 137-252.
\bibitem{Es} T. Estermann, On Kloosterman's sum, \emph{Mathematika} \textbf{8} (1961), 83-86.
\bibitem{FKM} E. Fouvry, E. Kowalski \& Ph. Michel, Algebraic trace functions over the primes. to appear in \emph{Duke Math. J.}
\bibitem{FM1} E. Fouvry \& Ph. Michel, Crible asymptotique et sommes de Kloosterman, Proc. Session in Analytic Number Theory and Diophantine Equations, Bonner Mathematische Schriften, Vol. 360, 2003.
\bibitem{FM2} E. Fouvry \& Ph. Michel, Sur le changement de signe des sommes de Kloosterman, \emph{Annals of Math.} \textbf{165} (2007), 675-715.
\bibitem{GMPY} D.A. Goldston, Y. Motohashi, J. Pintz \& C.Y. Y\i ld\i r\i m, Small gaps between primes exist, \emph{Proc. Japan Acad. Ser. A Math. Sci.} \textbf{82} (2006), 61-65.
\bibitem{Ka1} N.M. Katz, Sommes Exponentielles, Asterisque 79, Soci\'et\'e math\'ematique de France, 1980.
\bibitem{Ka2} N.M. Katz, Gauss sums, Kloosterman Sums, and Monodromy Groups, Annals of Mathematics Studies, Vol. 116, Princeton University Press, Princeton, NJ, 1988.
\bibitem{Ku}  N.V. Kuznetsov, The Petersson conjecture for cusp forms of weight zero and the Linnik conjecture. Sums of Kloosterman sums, \emph{Mat. Sb.} \textbf{111} (1980), 334-383.
\bibitem{Ma} K. Matom\"{a}ki, A note on signs of Kloosterman sums, \emph{Bull. Soc. Math. France} \textbf{139} (2011), 287-295.
\bibitem{Mi} Ph. Michel, Autour de la conjecture de Sato-Tate pour les sommes de Kloosterman. I, \emph{Invent. Math.} \textbf{121} (1995), 61-78.
\bibitem{Se} A. Selberg, Sieve methods, {\it Proc. Sympos. Pure Math.}, Vol. \textbf{XX}, 311-351, Amer. Math. Soc., Providence, R-I., 1971.
\bibitem{SF1} J. Sivak-Fischler, Crible \'{e}trange et sommes de Kloosterman, \emph{Acta Arith.} \textbf{128} (2007), 69-100.
\bibitem{SF2} J. Sivak-Fischler, Crible asymptotique et sommes de Kloosterman, \emph{Bull. Soc. Math. France} \textbf{137} (2009), 1-62.
\bibitem{We} A. Weil, On some exponential sums, \emph{Proc. Nat. Acad. Sci. U.S.A.} \textbf{34} (1948), 204-207.
\end{thebibliography}

\bigskip

\bigskip

\end{document}